\documentclass{IEEEtran}

\usepackage{graphicx}

\usepackage[cmex10]{amsmath}
\usepackage{amsfonts}
\usepackage{amssymb}
\usepackage{textcomp}

\usepackage{hyperref}

\usepackage{ulem}

\begin{document}

\title{Alternatives with stronger convergence than coordinate-descent iterative LMI algorithms}

\author{Emile~Simon and Vincent~Wertz
\thanks{E. Simon and V. Wertz are with the Pole of Applied Mathematics, ICTEAM Institute, Universit\'e Catholique de Louvain, 4 avenue Georges Lema\^{i}tre, 1348 Louvain-la-Neuve, Belgium (Corresponding author: simonemile@gmail.com, Tel/Fax. +32 10 47 88 10/ 21 80).}}
\maketitle

\begin{abstract}

In this note we aim at putting more emphasis on the fact that trying to solve non-convex optimization problems with coordinate-descent iterative linear matrix inequality algorithms leads to suboptimal solutions, and put forward other optimization methods better equipped to deal with such problems (having theoretical convergence guarantees and/or being more efficient in practice). This fact, already outlined at several places in the literature, still appears to be disregarded by a sizable part of the systems and control community. Thus, main elements on this issue and better optimization alternatives are presented and illustrated by means of an example.

\end{abstract}


\section{Main point}

Probably because of the success and large domination of convex and Linear Matrix Inequality (LMI) optimization in systems and control theory for the last two decades, it regularly happens that researchers try to solve non-convex optimization problems through the use of LMIs or convex subsets/approximations, as if using this path was adequate for all problems or the only acceptable possibility. The point of this brief note is to shed more light on this misunderstanding, and put forward other optimization approaches with strong convergence guarantees and very efficient in practice for non-convex problems.

We consider in particular the problems admitting a Bilinear Matrix Inequality (BMI) representation: involving product terms between variables, sometimes called `complicating variables'. A classical attempt to solve these problems is given by the following heuristic:

\begin{enumerate}
	\item Split the set of complicating variables in two subsets,\\ so that fixing either subset turn the BMI into an LMI.
	\item Choose/design initial values for one of the two subsets. 
  \item Fix the variables of the subset obtained at the previous step, and take the variables of the other subset as free optimization variables.
	\item Minimize the objective with LMI optimization.
	\item Repeat steps 3) and 4) iteratively, until the objective value reaches a given target or decreases less than a given accuracy.   
\end{enumerate}

This forms what can be described as a coordinate-descent iterative LMI algorithm (CDILMI), which is the most common kind of iterative LMI algorithms (ILMIs). The sequence of objectives generated by this type of scheme is only guaranteed to be monotonously non-increasing.
\subsection{Lack of convergence of CDILMIs}
The issue is that such heuristics will lead to `partial optimal' solutions \cite{WH76}: optimal in the `directions' of the two subsets of variables, but not in all directions. This means that algorithms of this type stop at solutions which are not locally optimal and are unable to follow directions that would improve the objective. The contingency where such an algorithm leads to a locally optimal solution should be very seldom: as outlined in a similar context in \cite{HM97}, it should almost never be the case. Note also that ILMIs typically stall once the sequence of solutions reaches a border of the BMI feasible set, while being nowhere close to a locally optimal solution.

Remark however that it is possible to develop ILMIs guaranteed to converge to locally optimal solutions: one such rare example is given in \cite{QGMD11}, where the proof of convergence does not rely on a coordinate-descent principle (other ingredients are used to ensure convergence). Anyhow, ILMIs will not maintain the guarantee of convergence towards globally optimal solutions nor the polynomial time complexity bound (unless $\mathcal{P=NP}$), both holding for convex/LMI problems. 

Instead of only trying to explore non-convex spaces with convex subsets or approximations, other alternatives that should be considered are local optimization methods (which moreover do not require additional optimization variables, much unlike convex approximations or relaxations).
\subsection{Alternative 1: Gradient-based methods}
The main difficulty is that many optimization problems in systems and control are non-smooth, so usual gradients do not exist everywhere. For instance, minimizing the $\mathcal{H}_\infty$ norm is a min-max objective function, i.e. locally Lipschitz, and the spectral abscissa is not even locally Lipschitz. To ensure convergence for these problems, it is necessary to consider non-smooth optimization methods.

In the current context, two methods clearly stand out: the open-source HIFOO \cite{BHLO06,BHLO06B} and commercial \texttt{hinfstruct} \cite{AN062} (\cite{AN061}) (both implemented under Matlab). We strongly advise the reader to consult these papers and related works. Both methods consider classical optimal reduced-order (so non-convex) controller designs for continuous-time Linear Time Invariant systems. While the range of problems covered by HIFOO is currently broader than that of \texttt{hinfstruct}, the latter method allows designing any control architecture made of conventional components as well as customized blocks.

More to the point, what matters are their underlying mechanisms to deal with non-smoothnesses. We give a simplified description as follows. 

HIFOO uses a conventional BFGS algorithm in its first phase, and then a random gradient sampling and bundle in the second phase. This converges based on the hypothesis that the function is differentiable almost everywhere, so there exist gradients near the current iterate which can be exploited. The mechanism behind \texttt{hinfstruct} is more elegant because it relies on extension sets of the Clarke sub-differential at each iteration. Thus, local information is completely available and is used to build a quadratic tangent model of the objective, efficiently optimized at the current iterate. It can be found in some places like \cite{A11} that \texttt{hinstruct} should in general be faster than HIFOO. 

Considering their theoretical convergence guarantees, ease of use, robustness and efficiency in practice, using either of these methods should present only advantages compared to ILMIs or any heuristics with weak convergence properties.
\subsection{Alternative 2: Derivative-free methods}
Another possible direction to solve optimization problems, which makes sense when gradient information is not available (at all, or accurately, or at an acceptable computational cost), is to consider Derivative-Free Optimization (DFO) methods (see the book \cite{CSV09} for a comprehensive presentation, or \cite[Chap. 6]{S12} for a summarized one). 

Apparently, the first work in systems and control where the idea of using a DFO method was thoroughly investigated is developed in the related papers \cite{AN04} and \cite{AN061}: The authors investigate the convergence of the multidirectional search (MDS) algorithm on non-smooth problems (spectral abscissa and $\mathcal{H}_\infty$ norm), outline the lack of convergence of this method for these problems, and propose additional non-smooth (first-order) steps which ensure convergence (this work later led to \texttt{hinfstruct} \cite{AN062}).

Next to \cite{AN04} and \cite{AN061}, the direction of using classical DFO methods was also visited in \cite{H06}, with the fundamental problem of static output feedback (SOF) stabilization (using three implementations available in \cite{MCT}). The benchmark results in \cite{H06} show that these methods are very successful for these (conjectured) $\mathcal{NP}$-hard problems. Motivated by both these papers and the availability of the proved convergent method HIFOO, we investigated further the performance of a DFO method on a broad benchmark of not only finding stabilizing SOF controllers but also minimizing the $\mathcal{H}_2$ and $\mathcal{H}_\infty$ norm of a performance channel \cite{S11}. It appeared that, despite only 0th order information was used by the DFO method, it performed reasonably well compared to HIFOO (sometimes better, both for the cpu times or objective values).

Still, as outlined in \cite{H06}, an explanation for the good performance of DFO methods remains to be identified. In our opinion, such reasons were partly outlined in \cite{AN061} when references are given to \cite{T97} and \cite{AD06}, where are developed important convergence results. It is worthwhile to note that these strong convergence guarantees are largely absent from the systems and control literature: we could not find any paper where a development relies on these guarantees. These proofs, given as follows, are however the key behind the convergence of DFO methods. Five results must be noted in particular: \cite{T97,DLT03,AD03,AD06,VC11}, summarized as follows.

\newpage

For smooth unconstrained problems, convergence to a stationary point is guaranteed by \cite{T97}, and the gradient norm is shown to be tied to the step size parameter \cite{DLT03},\cite[p. 123]{CSV09}. A hierarchical convergence analysis is proposed in \cite{AD03} for the non-smooth case: the authors provide convergence results for assumptions ranging from strict differentiability, regularity, Lipschitz continuity, lower semi-continuity and even general non-smooth functions. Their main result in the Lipschitz case is that the method generates a limit point where the Clarke directional derivatives are non-negative in a set of positive spanning directions. Later, \cite{AD06} generalizes the method and proposes the MADS algorithm for non-smooth constrained optimization problems, for which the set of directions with non-negative Clarke derivatives becomes asymptotically dense in $\mathbb{R}^n$. This convergence analysis is extended to the second-order in \cite{AA06} and for a class of discontinuous functions in \cite{VC11}. Note that since only a finite number of search directions are generated in practice, these asymptotical convergence guarantees may present some limitations (see \cite[pp. 111-114]{S12} and also \cite{A04}). 

From an user perspective, DFO methods are simple to use and implement and may quickly be tried as a first attempt to gain insights on an optimization problem or alternatively as last resort when other methods are not adequate. These methods, now supported by strong convergence properties, will yield better solutions than those obtained with heuristics not backed up by such convergence analysis. We recommend \cite[Chap. 1]{CSV09} for a more detailed description. We also refer to the thesis \cite{S12} and ref. therein for broader analysis and presentations on the current topic.

In summary, as long as an ILMI is not guaranteed to lead to locally optimal solutions, this kind of scheme might only be useful to find initial solutions (and not even necessarily then, because such solutions may be poorly located). Therefore, other optimization methods having better convergence properties and/or efficiency in practice should be used instead.

The main point of this note has been drawn. In the second (and last) section, we draw an illustration by means of an example with a CDILMI recently proposed.
\section{Example of CDILMI and results}
We illustrate the above-mentioned ideas with the algorithm of \cite{LLS10} which considers the design of reduced-order filters that must be positive (i.e. with all entries of the state-space matrices positive in the discrete-time case) and respecting a given performance level, with the performance chosen as the $\mathcal{H}_\infty$ norm of the filtering error system, in the context of discrete-time positive LTI systems.

It must be noted that in \cite{LLS10} the objective is to find a solution under a given performance level and not to minimize this objective function. In that sense, the focus there is not put on the convergence of the algorithm towards locally optimal solutions, but rather on several other contributions of the paper. In particular, a structure is put forward using the system augmentation approach which can be exploited to deal with constraints (here the positivity constraint) that could otherwise not be cast under a BMI and corresponding ILMI.
\vspace{0.55cm}

\subsection*{Numerical results and comments}
The details of the problem are given in \cite[Sec. IV]{LLS10}, and are omitted here for the sake of brevity\footnote{We only mention that the value of $b_3$ was erroneously written $=0.0128$ in \cite{LLS10}, instead of its correct value $=0.385$}. We directly copy the state-space matrices of the filter obtained in \cite{LLS10} hereunder:

\begin{center}
	$\hat{A} = 0.22819,\ \hat{B} = \left[0.00003\ 0.00003\right]$, 
	
	$\hat{C} = 0.14130,\ \hat{D} = \left[0.17889\ 0.34404\right]$.
\end{center}

This filter solution is not satisfying because its dynamical part (matrices $\hat{A}$, $\hat{B}$, $\hat{C}$) is canceled out by the very small entries of the $\hat{B}$ matrix. So this filter can be approximated by only its $\hat{D}$ matrix, without almost no impact ($0.003\%$) on the performance level (around 0.1417). What probably happened is that the CDILMI was blocked to this solution against a border of the feasible set (here the positivity constraint) which is a typical phenomenon of these type of algorithms.

This illustrates that more convergent local optimization methods should have been considered instead. For such LTI filter/controller design problems, we recommend in particular the methods HIFOO and \texttt{hinfstruct}, guaranteed to find locally optimal solutions and very efficient in practice. Note however that the current versions of these programs are implemented for continuous-time LTI systems and do not yet feature the possibility to add a constraint of an admissible range for the variables, such as the positivity constraint.  

Anyhow, using gradient information will not prove necessary to get excellent results on the considered problem.

For illustration purposes, we performed 100 optimizations from random initial solutions with several DFO methods\footnote{Full experimentations details are given on \cite[Subsec. 8.2.3]{S12}}. Almost all of the solutions had a lower performance level than 0.1415 and most of the solutions thus obtained (depending on the method used) had performance levels between 0.0447 and 0.0448, and never lower levels (thus 0.0447 is a candidate for a globally optimal level). This gives one illustration of the fact that general-purpose local optimization methods may largely outperform CDILMIs.

One of the obtained solutions is for instance the following:

\begin{center}
	$\hat{A} = 0.0561,\ \hat{B} = \left[0.2686\ 1.0749\right]$,
	
	$\hat{C} = 0.3094,\ \hat{D} = \left[0.1521\ 0.1089\right]$,
\end{center}
which is located inside the feasible set, i.e. not against the positivity constraint nor the limit of stability. We tried many local optimization methods from that solution and none could improve it, including a version of MADS which has a theoretical convergence guarantee applying to the current context of non-smooth locally Lipschitz problems \cite{AD03}.

One method in particular was between 90 and 100\% successful at reaching solutions with such performance level (depending on the accuracy used)\footnote{For convenience, the Matlab files reproducing these results are available on \url{www.mathworks.com/matlabcentral/fileexchange/33219}}: The implementation of the Nelder-Mead algorithm by N. Higham, available in \cite{MCT} and also used in \cite{H06}, which we restart at the last solution returned until no improvement is obtained to a given accuracy (this improves the method a lot, see the hyperlinked works for details). 
In practice, this method should lead to better solutions than ILMIs that are not supported by a strong convergence analysis. And when (sub)gradient information is readily available, we recommend using HIFOO or \texttt{hinfstruct} which come with solid convergence certificates of local optimality, or to develop a method based on the same mechanisms.

\small
\section*{Acknowledgments}
The authors gratefully acknowledge the reviewers and editor for their constructive comments and Charles Audet for its hand in the summary of DFO convergence results. Also acknowledged are Ping Li for comments on an earlier version and Didier Henrion for pointing out the reference \cite{HM97}. This research is supported by the Belgian Network DYSCO (Dynamical Systems, Control, and Optimization) funded by the Interuniversity Attraction Poles Programme of the Belgian State, Science Policy Office.
%


\end{document}